\input amstex
\magnification=\magstep1
\documentstyle{amsppt}
\advance\hoffset by 0.1 truein
\TagsOnRight
\hoffset=0.29in
\hsize=4.75in
\voffset=0in
\vsize=7.5in
\nopagenumbers
\document

\

\vskip1.0cm

\subheading{Linear systems and ramification points on 
reducible nodal curves}

\vskip0.5cm

\subheading\nofrills{\bf Eduardo Esteves}\footnote {Supported
by CNPq, Proc. 300004/95-8.}

\subheading\nofrills{\rm 
Instituto de Matem\'atica Pura e Aplicada, IMPA, Estrada 
Dona Castorina 110, 22460-320 Rio de Janeiro RJ, Brazil}

\vskip0.8cm

\baselineskip=12pt

\subheading\nofrills{\bf 1. Introduction}

\vskip0.4cm

\parindent=16pt

In \cite{\bf 1} Eisenbud and Harris developed a general theory in order to 
understand what happens to a linear system and its ramification points on a 
smooth curve when the curve degenerates to a curve $E$ of compact type. 
Eisenbud and Harris were able to obtain remarkable results from their 
theory, and we refer to \cite{\bf 1} for a description of some of these 
results and a partial list of the articles where these results are proved. 
In one of these articles Eisenbud and Harris asked: \sl What are the limits 
of Weierstrass points in families of curves degenerating to stable curves 
\underbar{not} of compact type?\rm \cite{\bf 2\rm , p. 499}. In the 
present note we hope to have found a satisfactory answer to the latter 
question (Theorem 7).

Actually, in the present note 
we deal with the more general situation of nodal curves $E$, not 
necessarily stable. We also deal with the degeneration of any linear 
system, not only the canonical system. Moreover, in contrast with the 
theory developed by Eisenbud and Harris, we do not need to blow up 
our degeneration family to swerve the degenerating ramification points 
away from the nodes of $E$. Indeed, we can assign the appropriate 
ramification weight to any node of $E$ (Theorem 7, item 2). Therefore, 
our note is a conceptual addition to the theory of 
limit linear systems even when $E$ is of compact type.

\vskip0.8cm

\subheading\nofrills{\bf 2. Degenerating linear systems}

\vskip0.4cm

\parindent=0pt

\bf Set-up: \rm 
Let $S$ be the spectrum of a discrete valuation ring $R$. Let $\pi$ be a 
parameter of $R$. Let $s$ 
(resp. $\eta$) denote the special (resp. generic) point of $S$, and let 
$\pi\in R$ be a parameter. 
Let $f\: C @>>> S$ be a flat, projective morphism. Suppose that the 
generic fibre $C(\eta)$ is a geometrically integral curve and the special 
fibre $C(s)$ is a nodal reduced curve. Assume that $C$ is a regular scheme. 
Let $C_1,\dots,C_t$ denote the irreducible components of $C(s)$. For 
each $C_i$ we let
$$
C_i^*:=C\setminus\bigcup_{j\neq i} C_j.
$$

\parindent=16pt

\vskip0.4cm

Since 
$C$ is regular, then $C_1,\dots,C_t$ are Cartier divisors, and any Cartier 
divisor on $C$ supported in $C(s)$ is a linear combination of $C_1,\dots,
C_t$. For every pair of integers $(i,j)$, with $1\leq i,j\leq t$, we let 
$\delta_{ij}$ denote the intersection number $C_i\cdot C_j$. It is clear 
that $\delta_{ij}$ is the number of points in $C_i\cap C_j$ if 
$i\neq j$. Since
$$
\Cal O_C\cong\Cal O_C(C_1+\dots+C_t),\tag{0.1}
$$
as $C_1+\dots+C_t$ is the Cartier divisor on $C$ cut out by $\pi$, then 
$$
\delta_{ii}=-\sum_{j\neq i} \delta_{ij}
$$
for every $i=1,\dots,t$. 
For each $i=1,\dots,t$, we will say that an effective Cartier divisor on $C$ 
supported in $C(s)$ is 
\sl $C_i$-free \rm if its support does not contain $C_i$.

Since $C$ is regular, for every invertible sheaf $L_{\eta}$ 
on $C(\eta)$ there is an invertible sheaf $\Cal L$ on $C$ such that 
$\Cal L(\eta)\cong L_{\eta}$. We call such an $\Cal L$ an 
\sl extension of $L_{\eta}$ to $C$. \rm 
It is easy to see that the sheaves $\Cal L\otimes
\Cal O_C(n_1C_1+\dots+n_tC_t)$ 
are all the extensions of $L_{\eta}$ to $C$. 

\vskip0.4cm

\parindent=0pt

\bf Set-up: \rm 
Fix an invertible sheaf $L_{\eta}$ on $C(\eta)$ of degree $d$, and a 
non-zero subvectorspace $V_{\eta}\subseteq H^0(C(\eta),L_{\eta})$ of 
rank $r+1$.

\vskip0.4cm

\parindent=16pt

If $\Cal L$ is an extension of 
$L_{\eta}$ to $C$, then we put:
$$
V_{\Cal L}:=V_{\eta}\cap H^0(C,\Cal L),
$$
where the above intersection is taken inside 
$H^0(C(\eta),L_{\eta})$. It is clear that $V_{\Cal L}$ is a 
free $R$-module of rank $r+1$ with $V_{\Cal L}(\eta)=V_{\eta}$. 
In addition, the induced homomorphism
$$
V_{\Cal L}(s) @>>> (f_*\Cal L)(s) @>>> H^0(C(s),\Cal L(s))
$$
is injective. To summarize, given an extension $\Cal L$ of $L_{\eta}$ to 
$C$, the linear system $(V_{\eta},L_{\eta})$ extends to the 
linear system $(V_{\Cal L},\Cal L)$ on $C$, 
whose restriction $(V_{\Cal L}(s),\Cal L(s))$ to $C(s)$ 
is a linear system. 

\vskip0.4cm

\parindent=0pt

\bf Definition: \rm We say that $(V_{\Cal L}(s),\Cal L(s))$ is a 
\sl limit linear system\rm .

\vskip0.4cm

\parindent=16pt

\proclaim{Theorem 1} For every irreducible component $C_i\subseteq C(s)$, 
there is a unique 
extension $\Cal L_i$ of $L_{\eta}$ to $C$ with the following 
properties:
\roster
\item the canonically induced homomorphism
$$
V_{\Cal L_i}(s) @>>> H^0(C_i,\left.\Cal L_i(s)\right|_{C_i})
$$
is injective;
\item if $\Cal I$ is an extension of $L_{\eta}$ to $C$ such that the 
induced homomorphism
$$
V_{\Cal I}(s) @>>> H^0(C_i,\left.\Cal I(s)\right|_{C_i})\tag{1.1}
$$
is injective, then there is an effective, 
$C_i$-free Cartier divisor $D$ on $C$ 
supported in $C(s)$ such that $\Cal I\cong\Cal L_i(D)$ and the induced 
homomorphism $V_{\Cal L_i} \hookrightarrow V_{\Cal I}$ is an isomorphism.
\endroster
\endproclaim

\demo{\bf Proof} We first show that there is an extension $\Cal I$ of 
$L_{\eta}$ to $C$ satisfying (1.1). In fact, choose any extension 
$\Cal J$ of $L_{\eta}$ to $C$. Let $n_1,\dots,n_t$ be integers such that
$$
\deg_{C_j} \Cal I(s)<0
$$
for every $j\neq i$, where $\Cal I:=\Cal J\otimes
\Cal O_C(n_1C_1+\dots+n_tC_t)$. 
Then
$$
V_{\Cal I}(s)\subseteq H^0(C(s),\Cal I(s))\subseteq \bigoplus_{j=1}^t 
H^0(C_j,\left.\Cal I(s)\right|_{C_j})=H^0(C_i,\left.\Cal I(s)\right|_{C_i}).
$$

Assume now that there are two extensions $\Cal I_1$ and $\Cal I_2$ of 
$L_{\eta}$ to $C$ such that
$$
V_{\Cal I_m}(s) @>>> H^0(C_i,\left.\Cal I_m(s)\right|_{C_i}) 
$$
is injective for $m=1,2$. We claim that there is an extension 
$\Cal N$ of $L_{\eta}$ to $C$ such that $\Cal I_m\cong\Cal N(D_m)$ for 
an effective, $C_i$-free Cartier divisor $D_m$ on $C$ supported in $C(s)$, and 
the induced homomorphism $V_{\Cal N} @>>> V_{\Cal I_m}$ is an isomorphism, 
for $m=1,2$. In fact, since both $\Cal I_1$ and $\Cal I_2$ are 
extensions of $L_{\eta}$, 
then there is a Cartier divisor $D$ on $C$ supported in $C(s)$ such that 
$\Cal I_1\cong\Cal I_2(D)$. It follows from (0.1) that we may assume that 
$D=D_1-D_2$, where $D_1$ and $D_2$ 
are disjoint, effective, $C_i$-free Cartier divisors on 
$C$ supported in $C(s)$. Let 
$$
\Cal M:=\Cal I_2(D_1) \text{\  \  and \  \  } 
\Cal N:=\Cal I_2(-D_2).
$$
It is clear that
$$
\Cal M\cong \Cal I_1(D_2) \text{\  \  and \  \  } 
\Cal N\cong \Cal I_1(-D_1).
$$
For $m=1,2$, the inclusion $\Cal I_m @>>> \Cal M$ induces the 
following commutative diagram:
$$
\CD
V_{\Cal M}(s) @>>> H^0(C_i,\left.\Cal M(s)\right|_{C_i})\\
@AAA @AAA\\
V_{\Cal I_m}(s) @>>> H^0(C_i,\left.\Cal I_m(s)\right|_{C_i}).
\endCD
$$
Since $D_1$ and $D_2$ are $C_i$-free, 
it follows that the right vertical homomorphism 
is an embedding for $m=1,2$. 
On the other hand, the bottom horizontal homomorphism is 
injective by assumption. It follows that the left vertical homomorphism 
is injective. Since $V_{\Cal M}$ and $V_{\Cal I_m}$ are free $R$-modules 
of same rank, then $V_{\Cal M}=V_{\Cal I_m}$ for $m=1,2$. On the 
other hand, since $D_1$ and $D_2$ are disjoint, then 
$V_{\Cal N}=V_{\Cal I_1}\cap V_{\Cal I_2}$ inside $V_{\Cal M}$. Therefore, 
$V_{\Cal N}=V_{\Cal I_m}$ for $m=1,2$. Our claim is proved. 
Note that $\Cal N\cong\Cal I_m$ if and only if $D_m=0$.

We finally show that there is a sheaf $\Cal L_i$ as in the statement of the 
theorem. Let $\Cal I_1$ be an extension of $L_{\eta}$ to $C$ such that 
the induced homomorphism
$$
V_{\Cal I_1}(s) @>>> H^0(C_i,\left.\Cal I_1(s)\right|_{C_i})
$$
is injective. If $\Cal I_1$ satisfies the second property in the statement 
of the theorem as well, then we are done: put $\Cal L_i:=\Cal I_1$. If not, 
by applying the reasoning in the above paragraph, there is an 
extension $\Cal I_2$ of $L_{\eta}$ and a non-zero, effective, 
$C_i$-free Cartier 
divisor $D_1$ on $C$ supported in $C(s)$ such that $\Cal I_2\cong
\Cal I_1(-D_1)$ and the induced homomorphism $V_{\Cal I_2} @>>> 
V_{\Cal I_1}$ is an isomorphism. It is clear that the induced 
homomorphism
$$
V_{\Cal I_2}(s) @>>> H^0(C_i,\left.\Cal I_2(s)\right|_{C_i})
$$
is injective. If $\Cal I_2$ satisfies the second 
property in the statement of the theorem as well, then we are done: put 
$\Cal L_i:=\Cal I_2$. If not, then procceed as before, thereby 
producing extensions $\Cal I_1,\Cal I_2,\dots,\Cal I_m$ on the $m$-th 
step such that $\Cal I_j\cong\Cal I_{j+1}(D_j)$ for a non-zero, effective, 
$C_i$-free Cartier divisor $D_j$ on $C$ supported in $C(s)$, and the induced 
homomorphism $V_{\Cal I_{j+1}} @>>> V_{\Cal I_j}$ is an isomorphism, for every 
$j<m$. In particular, we have that
$$
V_{\Cal I_1}=V_{\Cal I_m}
\subseteq H^0(C,\Cal I_1(-D_1-\dots-D_{m-1})).\tag{1.2}
$$
It follows from (1.2) that the above procedure cannot go on indefinitely. 
Thus there will be an $m\geq 1$ such that $\Cal L_i:=\Cal I_m$ is as in 
the statement of the theorem.

The uniqueness of $\Cal L_i$ is obvious from its properties. The proof 
of the theorem is complete.\qed
\enddemo

\vskip0.4cm

\parindent=0pt

\bf Definition: \rm 
We say that $\Cal L_i$ is the extension of $L_{\eta}$ \sl associated 
to $C_i$ \rm 
(and the subvectorspace $V_{\eta}\subseteq H^0(C(\eta),L_{\eta})$). We 
say that $(V_{\Cal L_i}(s),\Cal L_i(s))$ is the \sl limit linear 
system associated to $C_i$. \rm 

\parindent=16pt

\vskip0.4cm

\proclaim{Proposition 2} If $\Cal I$ is an extension of $L_{\eta}$ to $C$, 
then $\Cal I\cong\Cal L_i$ if and only if
\roster
\item the canonically induced homomorphism
$$
V_{\Cal I}(s) @>>> H^0(C_i,\left.\Cal I(s)\right|_{C_i})\tag{2.1}
$$
is injective;
\item for every irreducible component $C_j\subseteq C(s)$ with $j\neq i$, 
the canonically induced homomorphism
$$
V_{\Cal I}(s) @>>> H^0(C_j,\left.\Cal I(s)\right|_{C_j})\tag{2.2}
$$
is not identically zero.
\endroster
\endproclaim

\demo{\bf Proof} Suppose that $\Cal I\cong\Cal L_i$. Then (2.1) is 
injective. Suppose by contradiction that there is an 
irreducible component $C_j\subseteq C(s)$ with $j\neq i$ such that 
(2.2) is identically zero. It follows that $V_{\Cal I}=V_{\Cal J}$, where 
$\Cal J:=\Cal I(-C_j)$, contradicting the minimality property of $\Cal L_i$. 

Conversely, suppose that (2.1) is injective and (2.2) is not identically 
zero for every $j\neq i$. Since (2.1) is injective, then there is an 
effective, $C_i$-free Cartier divisor $D$ on $C$ supported in $C(s)$ such that 
$\Cal I\cong\Cal L_i(D)$ and the induced homomorphism $V_{\Cal L_i} @>>> 
V_{\Cal I}$ is an isomorphism. It follows that
$$
V_{\Cal I}\subseteq H^0(C,\Cal I(-D)),
$$
and hence every section of $V_{\Cal I}(s)$ is zero on $D$. Since 
(2.2) is not identically zero for every $j\neq i$ and $D$ is $C_i$-free, then 
$D=0$, and hence $\Cal I\cong\Cal L_i$.\qed
\enddemo

\vskip0.4cm

\parindent=0pt

\bf Remark 3. \rm Ziv Ran had also studied degenerations of linear 
systems in \cite{\bf 4}, where he had also obtained the linear system 
$(V_{\Cal L_i},\Cal L_i)$ of Theorem 1 for each component $C_i$. He called 
such system an ``effective state with focus $C_i$''.

\parindent=16pt

\vskip0.4cm 

\proclaim{Proposition 4} Fix $i,j$ with $i\neq j$. Let $l_{im}$, 
for $m\in\{1,\dots,t\}\setminus\{i\}$, be the unique integers such that
$$
\Cal L_i\cong\Cal L_j(\sum_{m\neq i} l_{im}C_m).
$$
Then $0\leq l_{im}\leq l_{ij}$ for every $m$.
\endproclaim

\demo{\bf Proof} Let
$$
\aligned 
E:=& \sum_{l_{im}>0} l_{im}C_m;\\
F:=& -\sum_{l_{im}<0} l_{im}C_m.
\endaligned
$$
We have that $\Cal L_i\cong\Cal L_j(E-F)$, where $E$ and $F$ are 
disjoint, effective, 
$C_i$-free Cartier divisors on $C$ supported in $C(s)$. Let
$$
\aligned 
\Cal M:=& \Cal L_j(E);\\
\Cal N:=& \Cal L_j(-F).
\endaligned
$$
It is clear that $\Cal M\cong\Cal L_i(F)$ and $\Cal N\cong\Cal L_i(-E)$. 
Since $F$ is $C_i$-free, the embedding $\Cal L_i \hookrightarrow \Cal M$ 
induces an isomorphism $V_{\Cal L_i}\cong V_{\Cal M}$. Since $E$ and 
$F$ are disjoint, then $V_{\Cal N}=V_{\Cal L_i}\cap V_{\Cal L_j}$ 
inside $V_{\Cal M}$. Thus $V_{\Cal N}\cong V_{\Cal L_j}$. So the 
induced homomorphism
$$
V_{\Cal L_j}(s) @>>> H^0(C(s),\Cal L_j(s))
$$
factors through the homomorphism
$$
H^0(C(s),\Cal N(s)) @>>> H^0(C(s),\Cal L_j(s))
$$
induced by the embedding $\Cal N \hookrightarrow \Cal L_j$. It follows from 
Proposition 2 that $F=0$. So $l_{im}\geq 0$ for every $m$.

On the other hand, we have that
$$
\Cal L_j\cong\Cal L_i(-\sum l_{im}C_m)\cong\Cal L_i(l_{ij}C_i+
\sum_{m\neq j}(l_{ij}-l_{im})C_m).
$$
Applying the result of the previous paragraph to the above situation, 
we obtain that $l_{ij}-l_{im}\geq 0$ for every $m\neq j$. The 
proof of the proposition is complete.\qed
\enddemo

\vskip0.4cm

\parindent=0pt

\bf Definition: \rm We say that $l_{ij}$ is the \sl connecting number 
\rm of $\Cal L_i$ and $\Cal L_j$.

\parindent=16pt

\vskip0.4cm

Note that $l_{ij}$ depends only on the specializations $\Cal L_i(s)$ and 
$\Cal L_j(s)$.

\vskip0.4cm

\proclaim{Corollary 5} Let $\Cal I$ be an extension of $L_{\eta}$ to $C$. 
If the canonically induced homomorphism
$$
V_{\Cal I}(s) @>>> H^0(C_m,\left.\Cal I(s)\right|_{C_m})
$$
is injective for $m=i,j$, then $\Cal L_i\cong\Cal L_j$.
\endproclaim

\demo{\bf Proof} It follows from Theorem 1 that there is an effective, 
$C_m$-free Cartier divisor $D_m$ on $C$ supported in $C(s)$ such 
that $\Cal I\cong\Cal L_m(D_m)$, for each $m=i,j$. Thus
$$
\Cal L_i\cong\Cal L_j(D_j-D_i).
$$
Since $D_i$ is $C_i$-free and $D_j$ is $C_j$-free, it 
follows easily from Proposition 4 that $D_i=D_j$. 
The proof is complete.\qed
\enddemo

\vskip0.4cm

\proclaim{\bf Proposition 6} 
Let $C_i,C_j\subseteq C(s)$ be two irreducible components intersecting at 
$p\in C_i\cap C_j$. For 
each $m=i,j$, let $\epsilon_0^m(p),\dots,\epsilon_r^m(p)$ be the increasing 
sequence of orders of vanishing at $p$ of the linear system 
$$
(V_{\Cal L_m}(s),\left.\Cal L_m(s)\right|_{C_m}).
$$
Then
$$
\epsilon^i_h(p)+\epsilon^j_{r-h}(p)\geq l_{ij}
$$
for every $h=0,\dots,r$.
\endproclaim

\demo{\bf Proof} The proof is analogous to the one given by 
Eisenbud and Harris in \cite{\bf 1\rm , Prop. 2.1, p. 348}.\qed
\enddemo

\vskip0.8cm

\subheading\nofrills{\bf 3. Ramification points}

\vskip0.4cm

\parindent=0pt

\bf Set-up: \rm Assume from now on that the characteristic of the 
residue field $k(s)$ is 0.

\parindent=16pt

\vskip0.4cm

Let $\omega$ be the canonical sheaf on $C$ relative to $S$. 
If $\Cal L$ is an extension of $L_{\eta}$ to $C$, then we can associate 
(cf. \cite{\bf 3}) to 
the inclusion $V_{\Cal L} \hookrightarrow H^0(C,\Cal L)$ a section
$$
s_{\Cal L}\in H^0(C,\Cal L^{\otimes r+1}\otimes\omega^{\otimes 
\binom{r+1}{2}}),
$$
called the \sl ramification section of the linear system 
$(V_{\Cal L},\Cal L)$\rm . 
Let $Z_{\Cal L}$ denote the zero scheme of $s_{\Cal L}$. The 
subscheme $Z_{\Cal L}\subseteq C$ is called the \sl ramification subscheme of 
the linear system $(V_{\Cal L},\Cal L)$ on $C$. \rm The intersection 
$Z_{\eta}:=Z_{\Cal L}\cap C(\eta)$ is the ramification subscheme of 
$(V_{\eta},L_{\eta})$. Thus $Z_{\Cal L}$ does not contain $C(\eta)$. For 
every $i=1,\dots,t$, let $n_i^{\Cal L}$ denote the multiplicity of $C_i$ in 
$Z_{\Cal L}$. It is clear that $s_{\Cal L}$ factors through a section
$$
s^*_{\Cal L}\in H^0(C,\Cal L^{\otimes r+1}\otimes\omega^{\otimes 
\binom{r+1}{2}}\otimes\Cal O_C(-n_1^{\Cal L}C_1-\dots-n_t^{\Cal L}C_t)).
$$
Moreover, 
the zero scheme $Z$ of $s^*_{\Cal L}$ is the unique relative Cartier divisor 
on $C$ over $S$ such that $Z_{\eta}=Z\cap C(\eta)$. Of course, 
$$
Z_{\Cal L}=Z+\sum_{i=1}^t n_i^{\Cal L}C_i.
$$
Note in particular that, if $V_{\Cal L}(s)\subseteq 
H^0(C_i,\left.\Cal L(s)\right|_{C_i})$ for a certain irreducible 
component $C_i\subseteq C(s)$, then $Z\cap C^*_i=Z_{\Cal L}\cap C^*_i$. 

\vskip0.4cm

\parindent=0pt

\bf Set-up: \rm Let $Z$ denote the relative Cartier divisor on $C$ over $S$ 
whose generic fibre $Z(\eta)$ is the ramification subscheme of 
$(V_{\eta},L_{\eta})$. We call $Z(s)$ the \sl limit ramification divisor\rm . 
(If $(V_{\eta},L_{\eta})$ is the canonical system, then we call $Z(s)$ the 
\sl limit Weierstrass divisor\rm .) 
For every $q\in Z(s)$, we let $w_q$ denote the 
weight of $q$ in $Z(s)$.

\parindent=16pt

\vskip0.4cm

\proclaim{Theorem 7} For each $i=1,\dots,t$, 
let $Z_i\subseteq C_i$ be the ramification subscheme of 
$$
(V_{\Cal L_i}(s),\left.\Cal L_i(s)\right|_{C_i}).
$$
Let $q\in C(s)$. For every irreducible component $C_i\subseteq C(s)$ 
containing $q$, we let $w^i_q$ denote the 
weight of $q$ in $Z_i$. Then:
\roster
\item if $q\in C^*_i$, then $w_q=w_q^i$.
\item if $q\in C_i\cap C_j$ for $i\neq j$, then
$$
w_q=w_q^i+w_q^j+(r-l_{ij})(r+1).
$$
\endroster
\endproclaim

\demo{\bf Proof} As we remarked before, we have that $Z\cap C^*_i=
Z_{\Cal L_i}\cap C^*_i$. On the other hand, it is clear that 
$Z_{\Cal L_i}(s)\cap C^*_i=Z_i\cap C^*_i$. Thus, if $q\in C^*_i$ then 
$w_q=w_q^i$.

Suppose now that $q\in C_i\cap C_j$ for $i\neq j$. For $m=i,j$, let 
$\omega_m$ denote the dualizing sheaf on $C_m$. We have a canonical 
embedding
$$
\omega_m @>>> \left.\omega(s)\right|_{C_m}
$$
of invertible sheaves whose cokernel has length 1 at $q$. Thus
$$
w_q^m=a_m-\binom{r+1}{2},\tag{7.1}
$$
where $a_m$ is the order of vanishing at $q$ of the restriction of 
$s_{\Cal L_m}$ to $C_m$, for $m=i,j$. 
For each $m=i,j$, let $b_m$ be the order of 
vanishing at $q$ of the restriction of $s^*_{\Cal L_m}$ to $C_m$. Since 
$Z$ is equal to the zero scheme of $s^*_{\Cal L_m}$ for $m=i,j$, then  
$s^*_{\Cal L_i}=s^*_{\Cal L_j}$ and 
$$
w_q=b_i+b_j.\tag{7.2}
$$
On the other hand, it is clear that 
$$
\aligned 
b_i=a_i-n_j^{\Cal L_i},\\
b_j=a_j-n_i^{\Cal L_j}.
\endaligned\tag{7.3}
$$
On one hand, since $s^*_{\Cal L_i}=s^*_{\Cal L_j}$, then 
$$
\Cal L_i^{\otimes r+1}\cong\Cal L_j^{\otimes r+1}\otimes
\Cal O_C(\sum_{h=1}^t(n_h^{\Cal L_i}-n_h^{\Cal L_j})C_h).\tag{7.4}
$$
On the other hand,
$$
\Cal L_i^{\otimes r+1}\cong\Cal L_j^{\otimes r+1}
((r+1)l_{ij}C_j+(r+1)E_{ij}),\tag{7.5}
$$
where $E_{ij}$ is a $C_i$-free and $C_j$-free effective Cartier divisor on $C$ 
with support in $C(s)$. Combining (7.4) and (7.5) we get that
$$
n^{\Cal L_i}_j+n^{\Cal L_j}_i=(r+1)l_{ij}.\tag{7.6}
$$
Combining (7.1), (7.2), (7.3) and (7.6), we have the equality in the 
statement (2) of the theorem. The proof is complete.\qed
\enddemo

\vskip0.4cm

\proclaim{\bf Corollary 8} For each $i=1,\dots,t$, let 
$d_i$ denote the degree of $\Cal L_i(s)$ on $C_i$. Then
$$
\sum_{i=1}^t d_i=d+\sum_{i<j}\delta_{ij}l_{ij}.
$$
\endproclaim

\demo{\bf Proof} The above formula follows from 
the Pl\"ucker formulas giving the degrees of $Z,Z_1,\dots,Z_t$ and 
the formula in item 2 of Theorem 7.\qed
\enddemo

\vskip0.4cm

\proclaim{\bf Corollary 9} Let $C_i,C_j\subseteq C(s)$ be 
irreducible components that intersect at a certain $q\in C_i\cap C_j$. For 
$m=i,j$, let $\epsilon_0^m(q),\dots,\epsilon_r^m(q)$ 
be the increasing sequence of orders of vanishing at $q$ of the linear 
system
$$
(V_{\Cal L_m}(s),\left. \Cal L_m(s)\right|_{C_m}).
$$
Then $q\not\in Z$ if and only if 
$$
\epsilon_h^i(q)+\epsilon_{r-h}^j(q)=l_{ij}
$$
for all $h=0,\dots,r$.
\endproclaim

\demo{\bf Proof} It follows from 
item 2 of Theorem 7 that $q\not\in Z$ if and only if 
$$
w_q^i+w_q^j = (l_{ij}-r)(r+1).
$$
On the other hand, since
$$
w_q^m=\sum_{h=0}^r (\epsilon_h^m(q)-h)
$$
for $m=i,j$, then
$$
w_q^i+w_q^j=\sum_{h=0}^r(\epsilon_h^i(q)+\epsilon_{r-h}^j(q))-r(r+1).
$$
Hence, it follows from Proposition 6 that
$$
w_q^i+w_q^j=(l_{ij}-r)(r+1)
$$
if and only if
$$
\epsilon_h^i(p)+\epsilon_{r-h}^j(p)=l_{ij}
$$
for all $h=0,\dots,r$. Combining the above two if-and-only-if 
statements we finish the proof.\qed
\enddemo

\vskip0.8cm

\subheading\nofrills{\bf 4. Cases}

\vskip0.4cm

\parindent=0pt

\bf Case 10. \rm (Curves of compact type) 
We do not have much control over the multidegrees 
of the limit linear systems. We know their ranges: it follows from 
Proposition 2 that
$$
0\leq \deg_{C_j} \Cal L_i(s)\leq d
$$
for all $i,j$. If 
$C(s)$ is of compact type, then Eisenbud and Harris \cite{\bf 1} 
developed a theory 
of limit linear series: for every irreducible component $C_i\subseteq C(s)$, 
they associated the unique limit linear system $(V_i,L_i)$ with 
$\deg_{C_j} L_i=d\delta_{ij}$ for all $j$, where $\delta_{ij}$ is 
the Kronecker symbol. 
(Of course, such choice may not be possible if $C(s)$ is not of compact 
type.) Eisenbud and Harris called the collection 
$\{(V_i,\left. L_i\right|_{C_i})| 1\leq i\leq t\}$ 
a crude limit series. It is also possible, 
in a way analogous to Theorem 7, to determine the limit ramification divisor 
$Z$ out of the crude limit series. (Though Eisenbud and Harris \cite{\bf 1} 
prefer to blow up 
the degeneration family at the nodes of the special fibre in order to 
swerve the degenerating ramification points away from these nodes.) 
Very seldom does the crude limit series coincide with 
the collection $\{(V_{\Cal L_i}(s),\left.\Cal L_i(s)\right|_{C_i})|
1\leq i\leq t\}$. 
In fact, if the above two collections coincide, 
then $(V_i,L_i)$ does not have base points at any of the nodes of 
intersection of two components.

\vskip0.4cm

\parindent=0pt

\bf Case 11. \rm (Planar curves) 
Let $E\subset\text{\bf P}^2$ be the planar curve of degree $d$ 
consisting of an irreducible nodal curve $Q$ of degree $d-1$ 
and one of its general secants $M$. Let 
$p_1,p_2,\dots,p_{d-1}$ be the points of $Q$ where $M$ intersects. 

\parindent=16pt

There are several ways we can view $E$ as a limit of smooth plane 
curves. In general, let $G(t)\in \text{\bf C}[x,y,z,t]$ be a polynomial 
in the variables $x,y,z,t$ that is homogeneous of degree $d$ in $x,y,z$. 
Let $C\subset \text{\bf P}^2\times\text{\bf A}^1$ be the 
zero scheme of $G$. Assume that $C$ is regular and flat over $\text{\bf A}^1$ 
in a neighbourhood of the fibre $C(0)$. Assume that 
the fibre $C(\lambda)$ is smooth for a general specialization 
$\lambda\in\text{\bf A}^1$. Finally, assume that $G(0)=F$.

We are concerned with computing the limit of the sets of 
inflection points on the curves $C(\lambda)$ as $\lambda$ tends to 0. 
(Of course, the inflection points of $C(\lambda)$ are the ramification 
points of the linear system of hyperplanes of $\text{\bf P}^2$ restricted to 
$C(\lambda)$.) It is clear that the restriction 
$|\left.\Cal O_{\text{\bf P}^2}(1)\right|_E|$ of the 
complete linear system of hyperplanes 
on $\text{\bf P}^2$ is a limit linear system. Moreover, 
it follows from the characterization given by Proposition 2 that 
$|\left.\Cal O_{\text{\bf P}^2}(1)\right|_E|$ is the limit linear system 
associated to the component $Q$ of $E$. This linear system has 
degree $d-1$ on $Q$ and 1 on the secant $M$. Since $Q$ has degree $d-1$, then 
the divisor of inflection points $Z_Q$ associated to $Q$ has degree 
$3(d^2-4d+3)$. To get the limit linear system 
associated with $M$, we twist $\left.\Cal O_{\text{\bf P}^2}(1)\right|_E$ 
by $\left.\Cal O_C(-M)\right|_E$. (Further twisting is not possible, since 
we would obtain an invertible sheaf with negative degree on $Q$.) Therefore, 
the limit linear system $(V,L)$ on $E$ 
associated to $M$ has degrees $\deg_Q L=0$ 
and $\deg_M L=d$. Restricting $(V,L)$ to 
$M$ we get a linear system of degree $d$ and rank 2 that yields a ramification 
divisor $Z_M$ of degree $3(d-2)$. 
Since the connecting number of $Q$ and 
$M$ is 1, then it follows from Theorem 7 that we have an extra 
weight contribution of 3 at each of the $p_1,p_2,\dots, p_{d-1}$. 
The upshot is that 
the limit of the ramification divisors on the smooth planar curves 
degenerating to 
$E$ is 
$$
Z=3p_1+3p_2+\dots+3p_{d-1}+Z_Q+Z_M.
$$

We note, however, that $Z_M$ depends on the particular degeneration to $E$. 
For an example, let $F:=x(y^2+x^2-z^2)$. For 
each pair $c=(c_1,c_2)\in\text{\bf A}^2$, let 
$G_c(t):=F+t(c_1y^3+c_2y^2z)$. Thus our degeneration depends on the 
parameter $c$. Computing the limit ramification divisor $Z_c$ on $E$, we 
get that
$$
Z_c=3p_1+3p_2+(0:y_1:z_1)+(0:y_2:z_2)+(0:y_3:z_3),
$$
where the $(y_i:z_i)$, for $i=1,2,3$, are the zeros of 
the polynomial 
$$
H:=c_1y^3+3c_2y^2z+3c_1yz^2+c_2z^3.
$$

\vskip0.4cm

\parindent=0pt

\bf Case 12. \rm (Curves with no singular ramification points) In case 
$C(s)$ is of compact type, it is possible that the limit ramification 
divisor $Z$ 
does not contain any singular points. (In fact, even if this is not 
the case, Eisenbud and Harris \cite{\bf 1} pointed out that it is possible to 
make all degenerating ramification points converge to smooth points of 
$C(s)$ replacing $C(s)$ by a semistably equivalent curve.) 
In general, though, singular points of the special fibre tend to attract 
ramification points. In fact, it is easy to show from Corollary 9 that, 
if the limit 
ramification divisor $Z$ contains only smooth points of $C(s)$, then:
\roster 
\item $C_i^*$ is smooth for every $i$;
\item $l_{ij}\geq r$ for all $i,j$;
\item $l_{ij}\delta_{ij}\leq d$ for all $i,j$.
\endroster 
In particular, if $r=d/2$ (the case of the canonical system), then 
$\delta_{ij}\leq 2$ for all $i,j$. 

\parindent=16pt

As a matter of fact, there should be more restrictions on 
$C(s)$, especially if we consider degenerations of the canonical system. 
For instance, if $C(s)$ is the union of only two components $C_1,C_2$, 
meeting at 
nodes $p_1,p_2$, then the limit Weierstrass divisor $Z$ has 
no singular points only if 
$\omega_i\cong\Cal O_{C_i}((g_i-1)(p_1+p_2))$ for $i=1,2$. (We denote by 
$\omega_i$ the canonical sheaf on $C_i$, and by $g_i$ the genus of $C_i$, 
for $i=1,2$.) It would be interesting to have a characterization of the 
nodal curves whose canonical limit Weierstrass divisors do not include 
singular points.

\enddocument